\documentclass[runningheads]{llncs}
\usepackage{amsmath}
\usepackage{amssymb}
\usepackage[applemac]{inputenc}
\usepackage[colorlinks]{hyperref}
\usepackage{url}
\usepackage{graphicx,epstopdf}
\usepackage{subcaption}
 \usepackage{epstopdf}
\usepackage[applemac]{inputenc}
\usepackage[colorlinks]{hyperref}
\usepackage{url}
\usepackage{subcaption}
\usepackage{hyperref}
\usepackage[T1]{fontenc}
\usepackage{graphicx}
\begin{document}
\title{On canonical parameterizations of $2D$-shapes\thanks{Supported by FWF grant I 5015-N, Institut CNRS Pauli, Vienna, Austria, and University of Lille, France}}

%
%
\author{Alice~Barbora~Tumpach\inst{2, 3}\orcidID{0000-0002-7771-6758}}
\authorrunning{Tumpach}
%
\institute{
Institut CNRS Pauli, Oskar-Morgenstern-Platz 1, 1090 Vienna, Austria \and
University of Lille, Cit\'e scientifique, 59650 Villeneuve d'Ascq, France
\email{alice-barbora.tumpach@univ-lille.fr}\\
\url{http://math.univ-lille1.fr/~tumpach/Site/home.html}}

\maketitle              
\begin{abstract}
This paper is devoted to the study of unparameterized simple curves in the plane. 
We propose diverse canonical parameterizations of a 2D-curve. For instance, the arc-length parameterization is canonical, but we consider other natural parameterizations like the parameterization proportionnal to the curvature of the curve. Both aforementionned parameterizations are very natural and correspond to a natural physical movement : the arc-length parameterization corresponds to travelling along the curve at constant speed, whereas parameterization proportionnal to curvature corresponds to a constant-speed moving frame. Since the curvature function of a curve is a geometric invariant of the unparameterized curve, a parameterization using the curvature function is a canonical parameterization.  The main idea is that to any physically meaningful stricktly increasing function is associated a natural parameterization of 2D-curves, which gives an optimal sampling, and which can be used to compare unparameterized curves in a efficient and pertinent way. An application to point correspondance in medical imaging is given.

\end{abstract}
\section{Introduction}

Curves in $\mathbb{R}^2$  appear in many applications: in shape recognition as outline of an object, in radar detection as the signature of a signal, as trajectories of cars etc... There are two main features of the curve : the route and the speed profil. In this paper, we are only interested in the route drawn by the curve and we will called it the unparameterized curve. An unparameterized curve can be parameterized in multiple ways, and the choosen parameterization selects the speed at which the curve is traversed. Hence a curve can be travelled with many different speed profils, like a car can travel with different speeds (not necessarily constant) along a given road. The choice of a speed profil is called a parameterization of the curve. It may be physically meaningful or not. For instance, depending on applications, there may not be any relevant parameterization of the contour of the statue of Liberty depicted in Fig.~\ref{Liberty3}. In this paper, we propose various very natural parameterization of 2D-curves.
They are based on the curvature, which together with the arc-length measure form a  complete set of geometric invariants or descriptors of the unparameterized curves. 

\section{Different parameterizations of 2D-shapes}\label{interpolate}

\subsection{Arc-length parameterization and Signed curvature}

 By 2D-shape, we mean the shape drawn by a parameterized curve in the plane. It is the ordered set of points visited by the curve. The shapes of two curves are identical if one can reparameterize one curve into the other (using a continuous increasing function). Any rectifiable planar curve admits a canonical parameterization, its \textit{arc-length parameterization}, which draws the same shape, but with constant unit speed. The set of 2D-shapes can be therefore identified with the set of arc-length parameterized curves, which is not a vector subspace, but rather an infinite-dimensional submanifold of the space of parameterized curves (see \cite{TumPres}). 
 
\vspace{-1cm} 
 \begin{figure}[!ht]
 		\centering
		\includegraphics[width = 7cm]{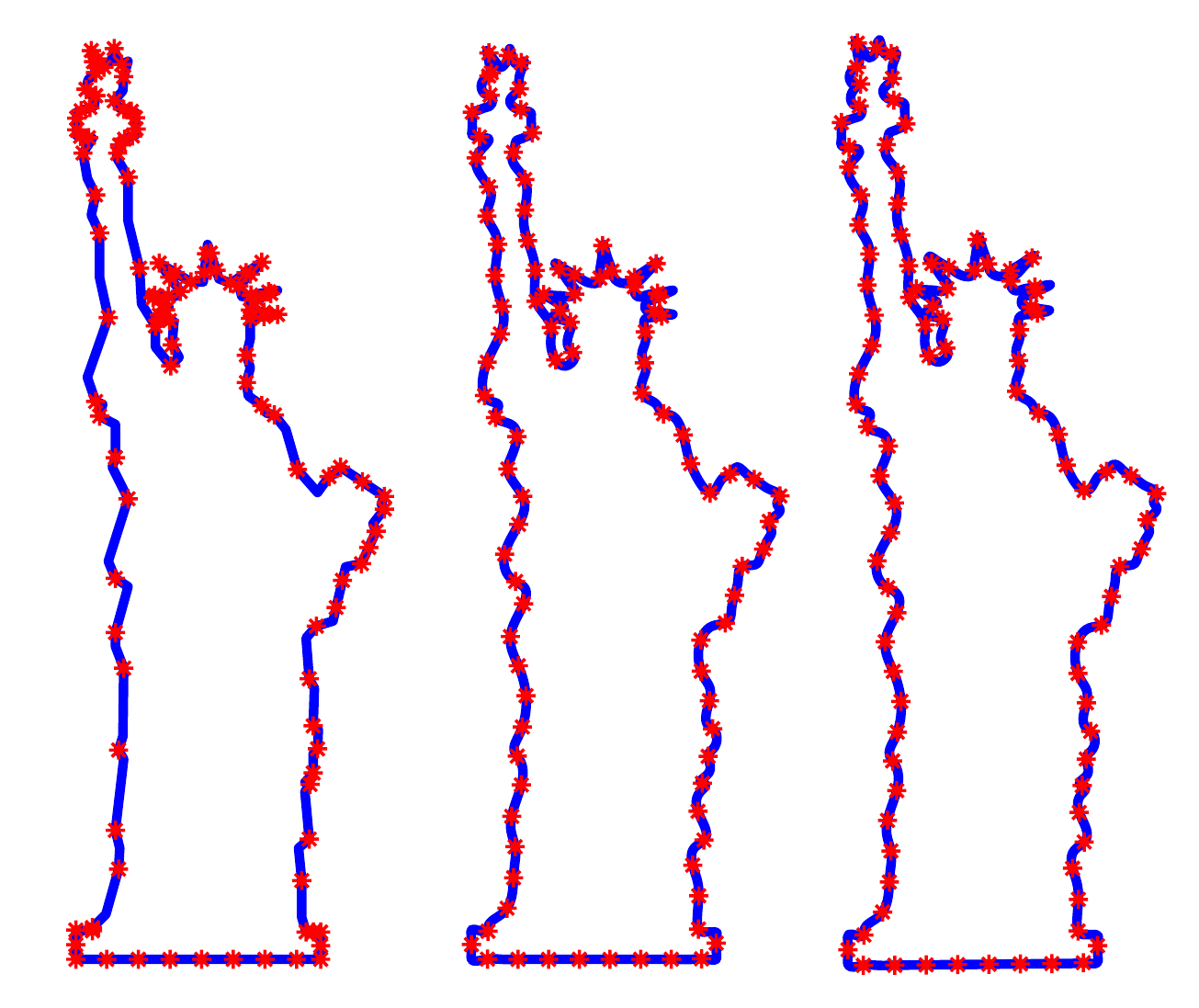}
 		\caption{\scriptsize{The statue of Liberty (left), a uniform resampling using Matlab function spline (middle), a reconstruction of the statue using its discrete curvature (right).}
		}
 		\label{Liberty3}		
\end{figure}

It may be difficult to compute an explicit formula of the arc-length parameterization of a given rectifiable curve. Fortunately, when working with a computer, one do not need it. One neither need a concrete parameterization of the curve to depict it, a sample of points on the curve suffises. To draw the statue of Liberty as in Fig.~\ref{Liberty3} left, one just need a finite ordered set of points (the red stars). The discrete version of an arc-length parameterized curve is a uniformly sampled curve, i.e. an ordered set of equally distant points (for the euclidean metric). Resampling a curve uniformly is immediate using some appropriate interpolation function like the matlab function \textit{spline} (the second picture in Fig.~\ref{Liberty3} shows a uniform resampling of the statue of liberty).

\begin{figure}[!ht]
 		\centering
		
 		\includegraphics[width = 8cm]{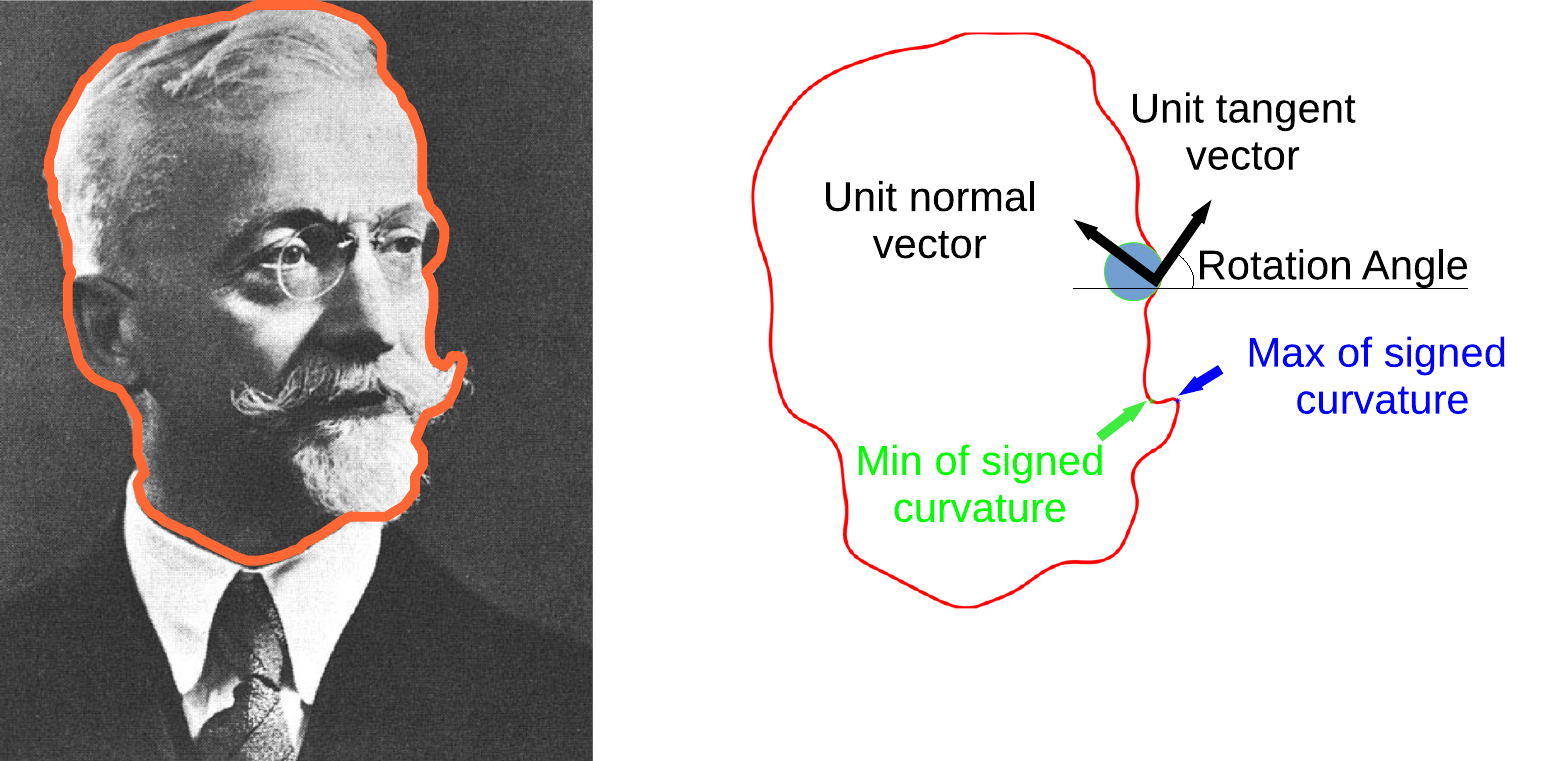}
		
 		\caption{\scriptsize Elie Cartan and the moving frame associated to the contour of  his head.}
 		\label{Cartans_head}	
\end{figure}

\begin{figure}[!ht]
 		\centering
		
		\includegraphics[width = 0.8\textwidth]{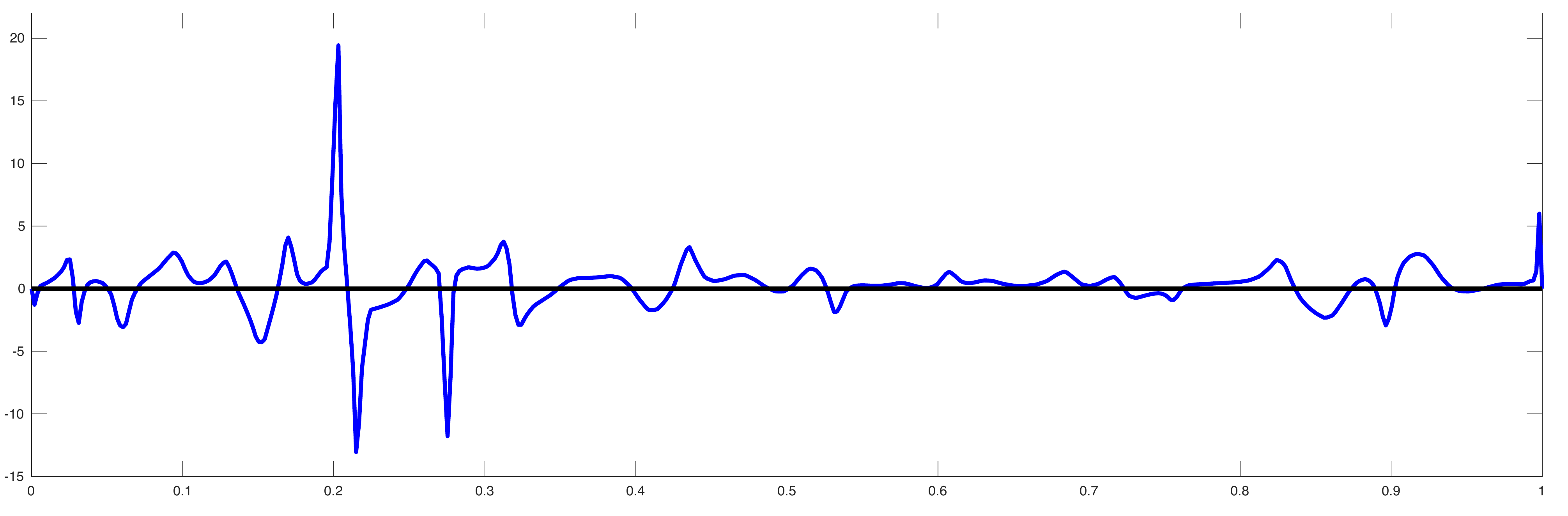}\\
		\includegraphics[width = 0.8\textwidth]{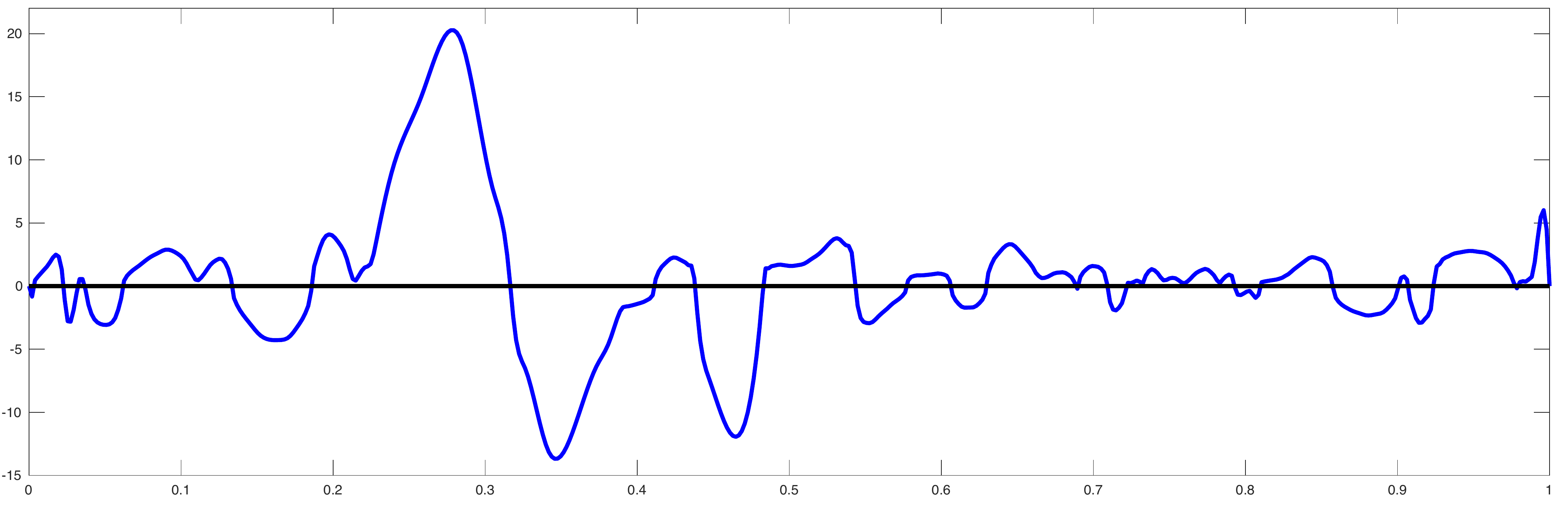}\\
		\includegraphics[width = 0.8\textwidth]{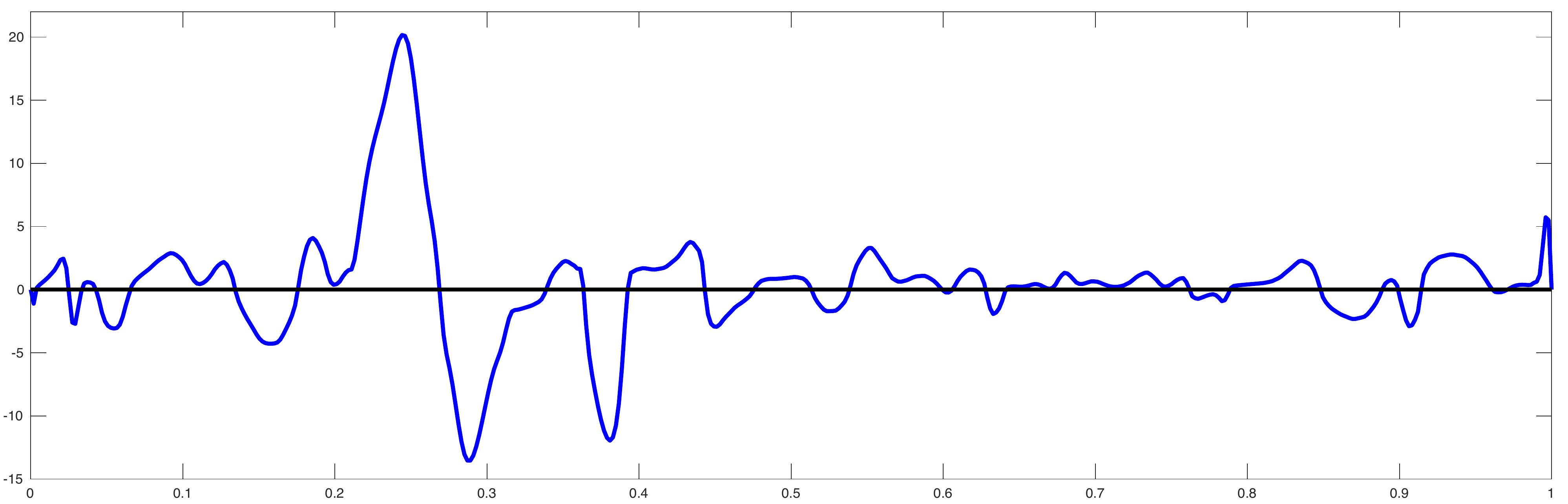}
		 		\caption{\scriptsize Signed curvature of Elie Cartan's head for the  parameterization proportional to arc-length (first line), proportional to the curvature-length (second line), and proportional to the curvarc length (third line).}
 		\label{Cartan_curvature}		
\end{figure}

Consider the set of 2D  simple  closed curves, such as the contour of Elie Cartan's head in Fig.~\ref{Cartans_head}. After the choice of a starting point and a direction, there is a unique way to travel the curve at unit speed. In Fig.~\ref{Cartans_head}, we have drawn the velocity vector near the glasses of  Elie Cartan, as well as the unit normal vector which is obtained from the unit tangent vector by a rotation of $+\frac{\pi}{2}$. These two vectors form an orthonornal basis, i.e. an element (modulo the choice of a basis of $\mathbb{R}^2$) of the Lie group $\textrm{SO}(2)$, which is characterized by a rotation angle. The rate of variation of this rotation angle is called the signed curvature of the curve. For instance, when moving along the external outline of the glasses, this curvature equals the inverse of the radius of the glasses. We have depicted the curvature function $\kappa$ of Elie Cartan's head in Fig.~\ref{Cartan_curvature}, first line, when the parameter $s\in[0; 1]$ on the horizontal axis is proportional to arc-length, and such that the entire contour of Elie Cartan's head is travelled when the parameter reaches 1. Its corresponds to a uniform sampling of the contour. The curvature function is also depicted when parameterized by two other canonical parameters, namely by the curvature-length parameter (second line) and the curvarc-length parameter (third line).

\vspace{-0.5cm}
\vspace{-0.2cm}

A discrete version of an arc-length parameterized curve is an equilateral polygon. To draw an equilateral polygon, one just need to know the length of the edges, the position of the first edge, and the angles between two successive edges. The sequence of turning-angles is the discrete version of the curvature and  defines a equilateral polygon modulo rotation and translation. In Fig.~\ref{Liberty3} , right, we have reconstructed the statue of Liberty using the discrete curvature.

In order to interpolate between two parameterized curves, it is easier when the domains of the parameter coincide. For this reason we will always consider curves parameterized with a parameter in $[0; 1]$. A natural parameterization is then the parameterization proportional to arc-length. It is obtain from the parameterization by arc-length by dividing the arc-length parameter by the length of the curve $L$. The corresponding curvature function is  also defined  on $[0; 1]$ and is obtained from the curvature function parameterized by arc-length by compressing the $x$-axis by a factor $L$. To recover the initial curve from the curvature function associated to the parameter $s\in[0; 1]$ proportional to arc-length, one only need to know the length of the curve.

\subsection{Parameterization proportional to curvature-length}
\begin{figure}[!ht]
 		\centering
		\includegraphics[width = 4cm]{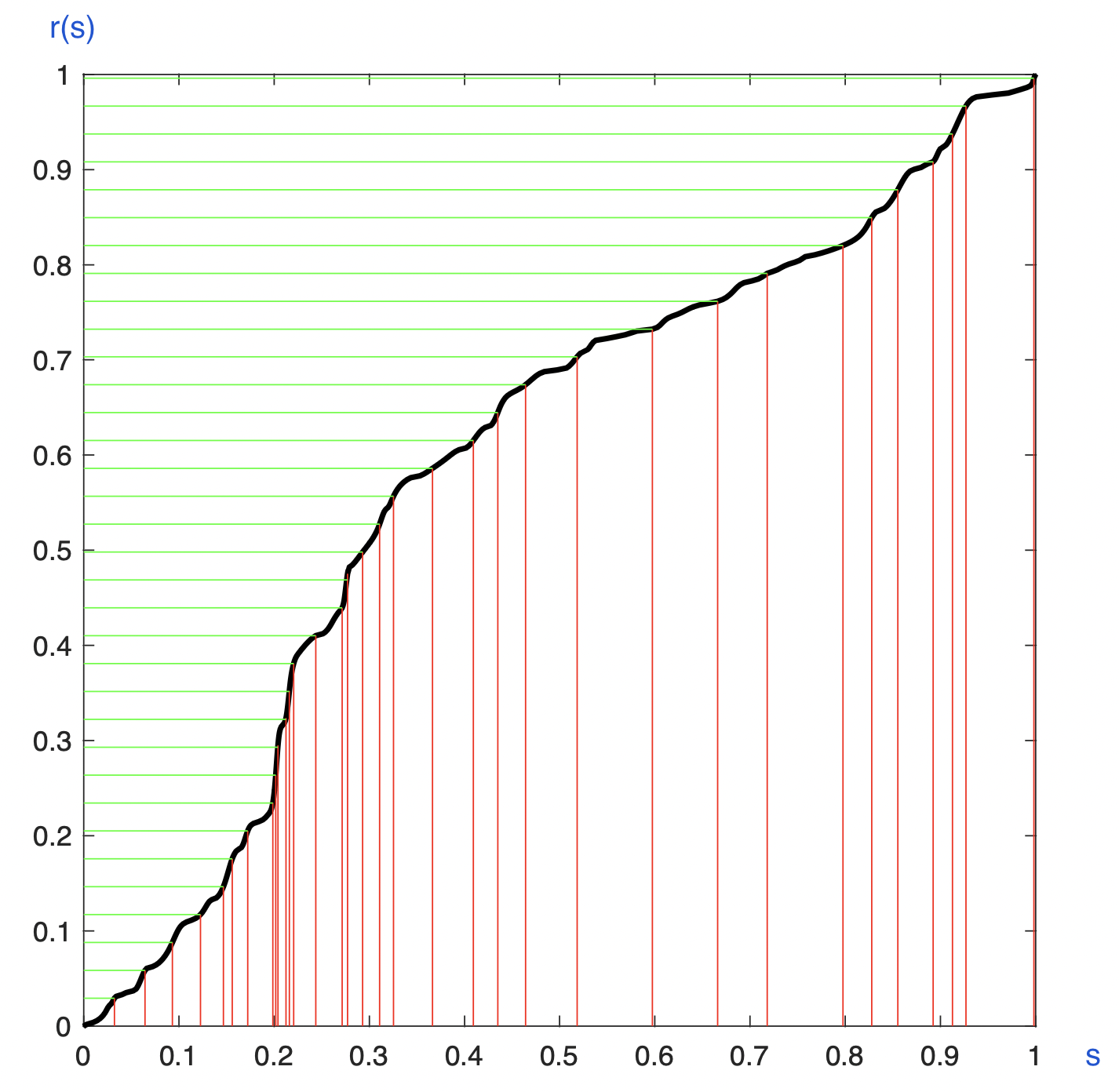} \hspace{0.5cm}
		\includegraphics[width = 4cm]{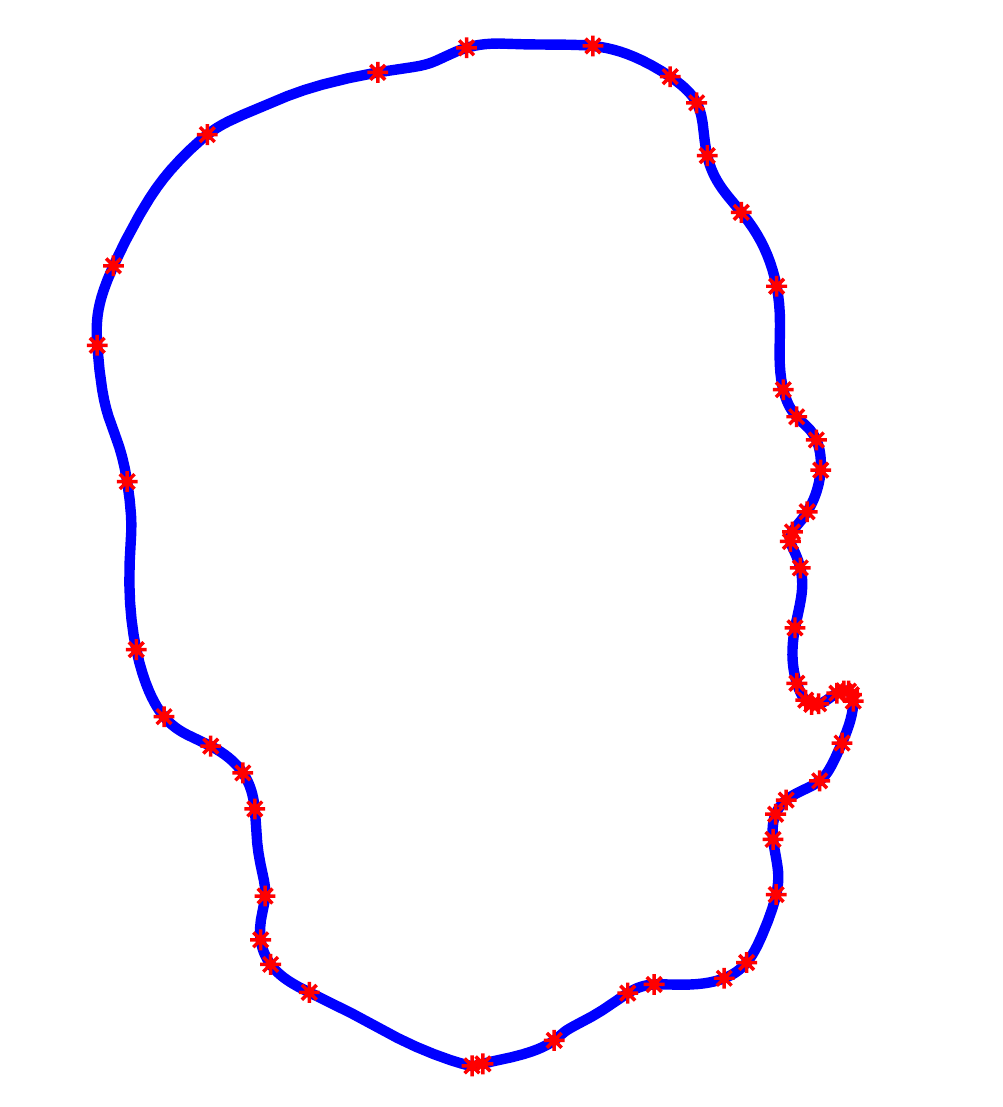}\\
		\vspace{-0.5cm}
		\includegraphics[width = 4cm]{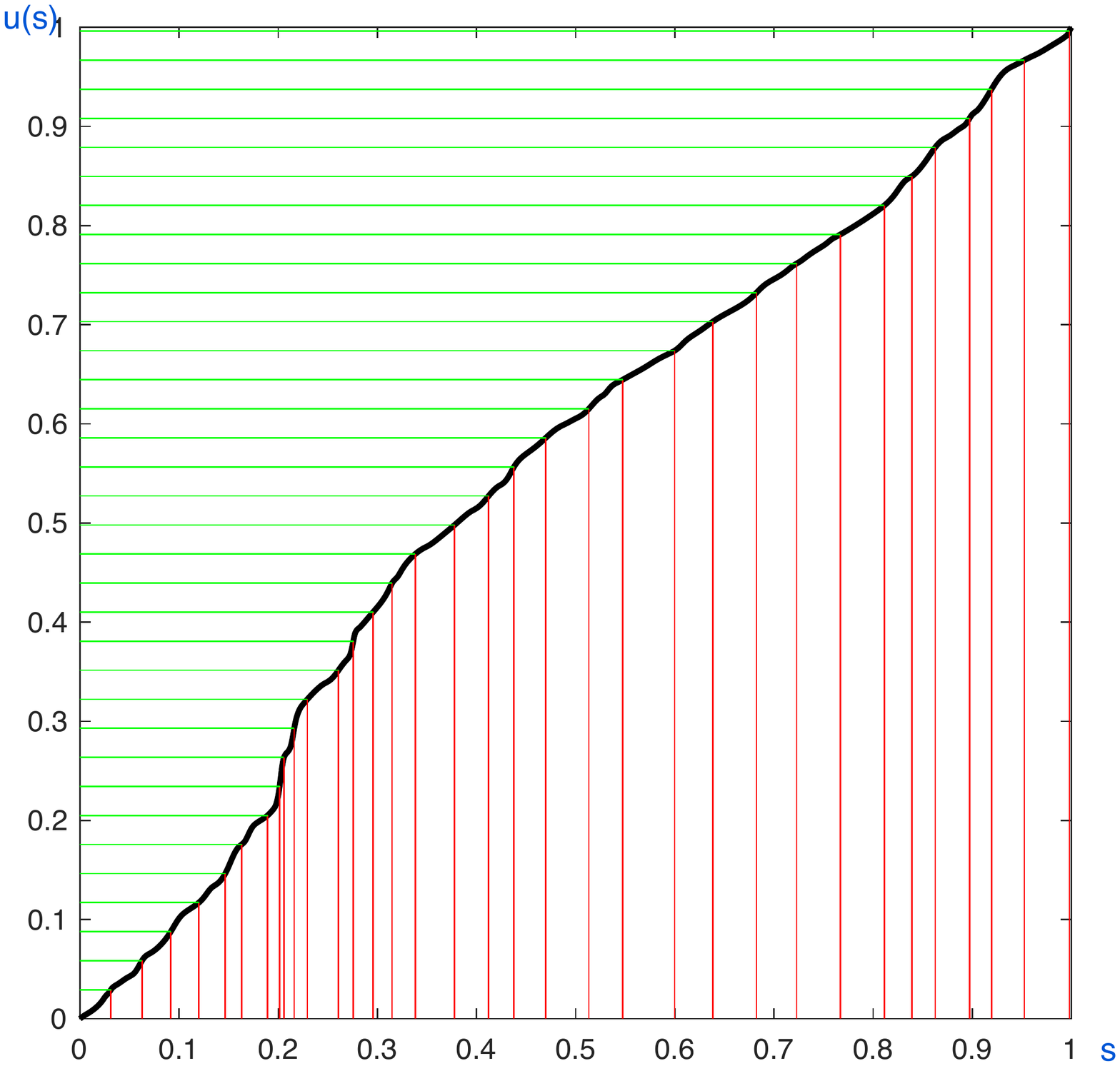}\hspace{0.5cm}
		\includegraphics[width = 4cm]{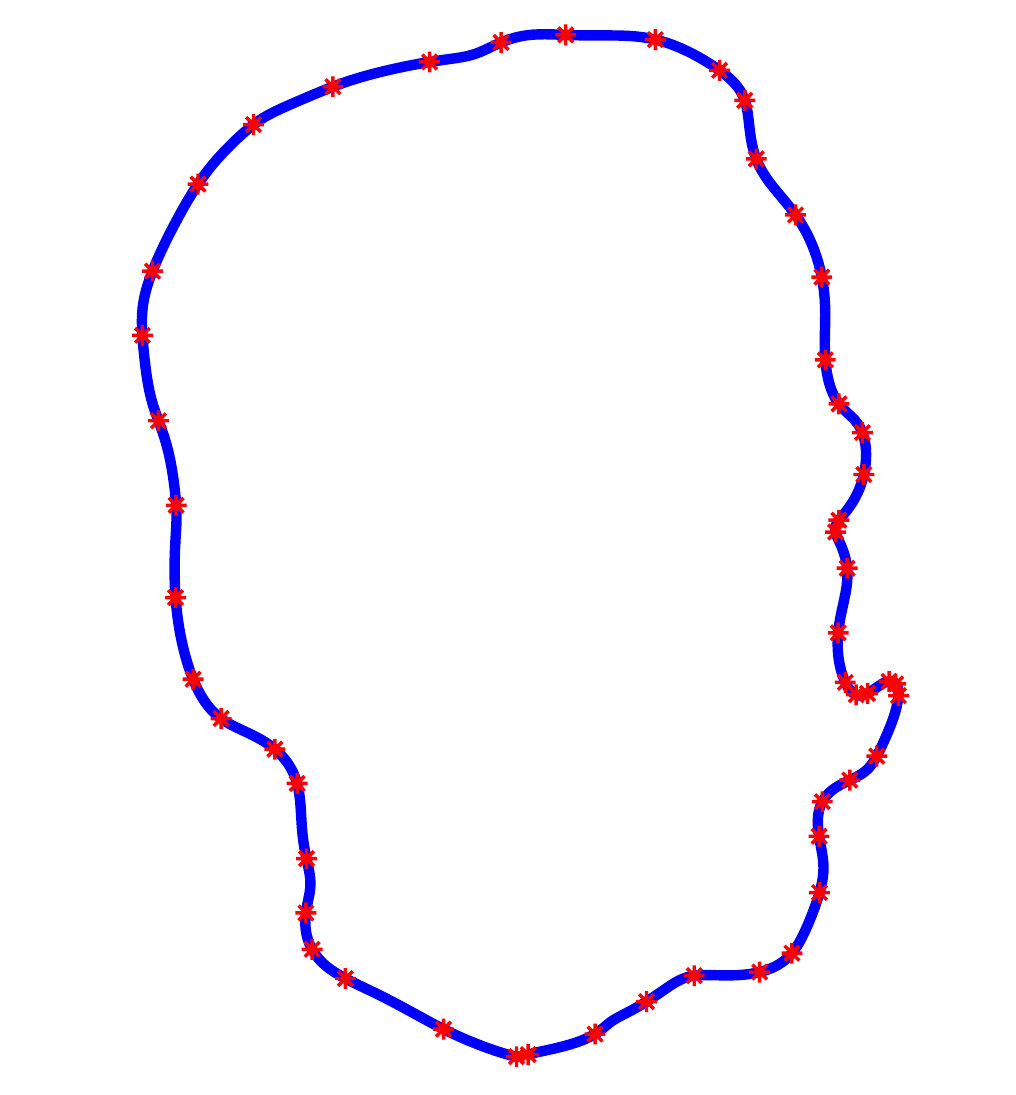}
 		\caption{\scriptsize First line : Integral of the (renormalized) absolute value of the curvature  (left), and corresponding resampling of Elie's Cartan head (right).
		Second line : Integral of the (renormalized) curvarc length  (left), and corresponding resampling of Elie's Cartan head (right).}
 		\label{Cartan_int_curvature}		
\end{figure}

In the same spirit as the scale space of T. Lindeberg (\cite{Lindeberg}), and the curvature scale space of Mackworth and Farcin Mokhtarian (\cite{Mokhtarian}), we now define another very natural parameterization space of 2D curves. Its relies on the fact that the integral of the absolute value of the curvature $\kappa$ is an increasing function on the interval $[0; 1]$, 
stricktly increasing when there are no flat pieces. In that case the function 
\begin{equation}\label{r}
r(s) = \frac{\int_0^s |\kappa(s)| ds}{\int_0^1 |\kappa(s)| ds}
\end{equation}
(where $\kappa$ denotes the curvature of the curve) belongs to the group of orientation preserving diffeomorphisms of the parameter space $[0; 1]$, denoted by $\textrm{Diff}^+([0; 1])$. Note that its inverse $s(r)$ can be computed graphically using the fact that its graph is the symmetric of the graph of $r(s)$ with respect to $y = x$. 
The contour of Elie Cartan's head can  be reparameterized using the parameter $r \in [0; 1]$ instead of the parameter $s\in [0; 1]$. In Fig.~\ref{Cartan_int_curvature} upper left, we have depicted the graph of the function $s\mapsto r(s)$. A uniform sampling with respect to the parameter $r$ is obtain by uniformly sampling the vertical-axis (this is materialized by the green equidistributed horizontal lines) and resampling Elie Cartan's head at the sequence of values of the $s$-parameter given by the abscissa of the corresponding points on the graph of $r$ (where  a green line hits the graph of $r$ a red vertical line materializes the corresponding abscissa). One sees that this reparameterization naturally increases the number of points where the 2D contour is the most curved, and decreases the number of points on nearly flat pieces of the contour. For a given number of points, it gives an optimal way to store the information contained in the contour.
The quantity
\begin{equation}\label{curvature_length}
C = L \int_0^1 |\kappa(s)| ds,
\end{equation}
where $s\in[0; 1]$ is proportional to arc-length, is called the {\it total curvature-length} of the curve. It is the length of the curve drawn in $\textrm{SO}(2)$ by the moving frame associated with the arc-length parameterized curve. 
For this reason we call this parameterization  the {\it parameterization proportional to curvature-length}. In the right picture of Fig.~\ref{Cartan_int_curvature}, we show the corresponding resampling of the contour of Elie Cartan's head. 

This resampling can naturally be adapted in the case of flat pieces resulting in a sampling where there is no points between two points on the curve joint by a straight line. In the left picture of Fig.~\ref{Liberty5}, we have depicted a sampling of the  statue of Liberty proportional to curvature-length. Note that there are no points on the base of the statue. The corresponding parameterization has the advantage of concentrating on the pieces of the contour that are very complex, i.e. where there is a lot of curvature, and not distributing points on the flat pieces which are easy to reconstruct (connecting two points by a straight line is easy, but drawing the moustache of Elie Cartan is harder and needs more information).


The drawback of using the parameterization proportional to curvature-length is that one can not reconstruct the flat pieces of a shape without knowing their lengths (remember that the parameterization proportional to curvature-length put no point at all on flat pieces). For this reason we propose a parameterization intermediate between arc-length parameterization and curvature-length parameterization. We call it {\it curvarc-length parameterization}.

\subsection{Curvarc-length parameterization}

\begin{figure}[!ht]
 		\centering
		\includegraphics[width = 0.9\textwidth]{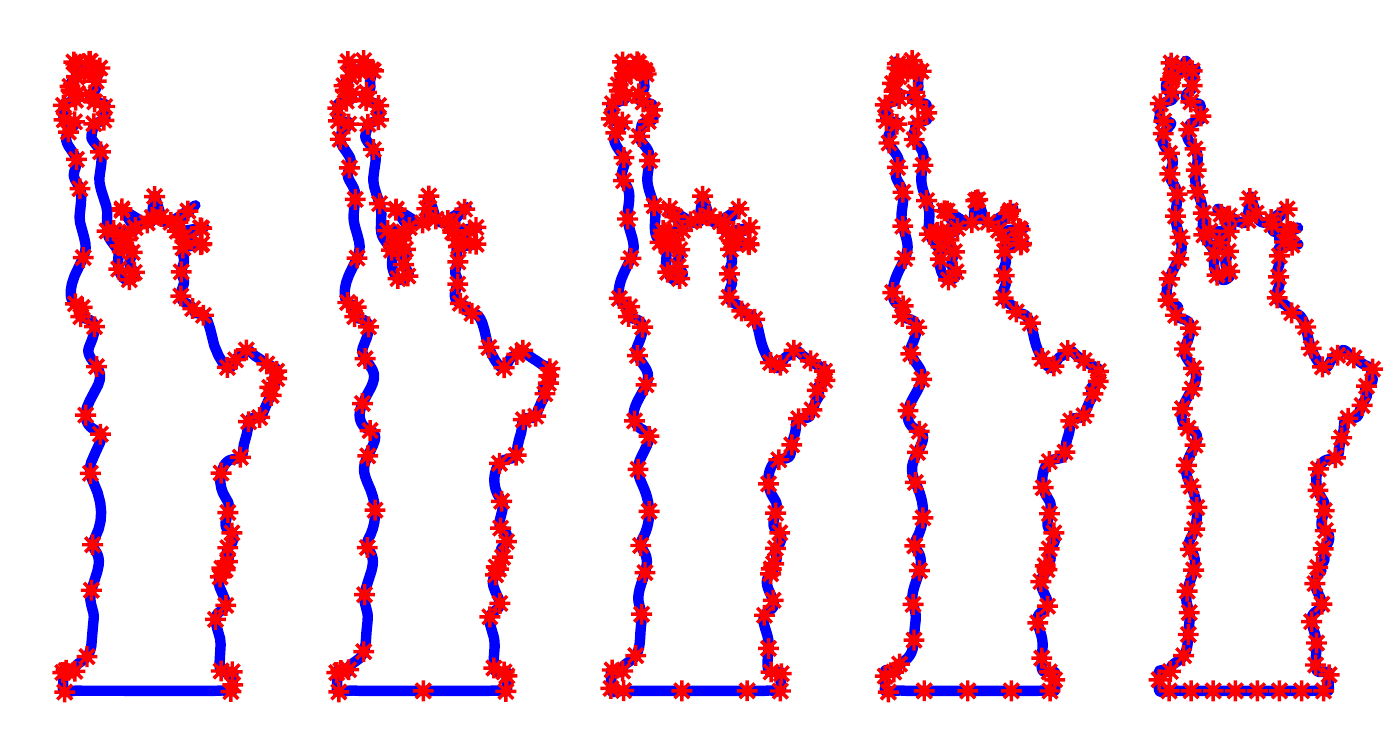}
 		\caption{\scriptsize 
		Resampling of the statue of Liberty proportional to the intergral of $\lambda$ + curvature, for (from left to right)
		$\lambda = 0; \lambda = 0.3; \lambda = 1; \lambda = 2; \lambda = 100$.		}
 		\label{Liberty5}		
\end{figure}

In order to define the curvarc-length parameterization, we consider the triple $(P(s), \vec{v}(s), \vec{n}(s))$, where $P(s)$ is the point of the shape parameterized proportionally to arc-length with $s\in[0; 1]$, $\vec{v}(s)$ and $\vec{n}(s)$ the corresponding unit tangent vector and unit normal vector respectively. It defines an element of the group of rigid motions of $\mathbb{R}^2$, called the special Euclidean group and denoted by $\textrm{SE}(2) := \mathbb{R}^2\rtimes \textrm{SO}(2)$. The point $P(s)$ corresponds to the translation part of the rigid motion, it is the vector of translation needed to move the origin to the point of the curve corresponding to the parameter value $s$. The moving frame $O(s)$ defined by $\vec{v}(s)$ and $\vec{n}(s)$ is the rotation part of the rigid motion. 
One has the following equations~:
\begin{equation}\label{triedre}
\frac{dP}{ds}  = L \vec{v}(s) \quad\textrm{and} \quad O(s)^{-1}\frac{d}{ds}O(s) = \left(\begin{smallmatrix} 0 & -\kappa(s)\\ \kappa(s) & 0\end{smallmatrix}\right),
\end{equation}
where $L$ is the length of the curve. 
Endow $\textrm{SE}(2) := \mathbb{R}^2\rtimes \textrm{SO}(2)$ with the structure of a  Riemannian manifold, product of the plane and the Lie group $\textrm{SO}(2) \simeq \mathbb{S}^1$. Than the norm of the  tangent vector to the curve $s\mapsto (P(s), \vec{v}(s), \vec{n}(s))$ is $L + |\kappa(s)|$. Therefore the length of the $\textrm{SE}(2)$-valued curve is $L + \int_0^1 |\kappa(s)| ds = L + \frac{C}{L}$. We call it the total curvarc-length. It follows that the following function
\begin{equation}
u(s) = \frac{\int_0^s (L + |\kappa(s)|) ds}{\int_0^1 (L + |\kappa(s)|) ds}
\end{equation}
defines a reparameterization of $[0; 1]$.  
More generally, one can use the following canonical parameter to reparameterize a curve in a canonical way:
\begin{equation}\label{ulambda}
u_\lambda(s) = \frac{\int_0^s L\lambda + |\kappa(s)|) ds}{\int_0^1 L\lambda + |\kappa(s)|) ds},
\end{equation}
where $s$ is the arc-length parameter.
In Fig.~\ref{Liberty5}  we show the resulting sampling of the Statue of Liberty for different values of $\lambda$. Note that for $\lambda = 0$, one recovers the curvature-length parameterization, for $\lambda = 1$ one obtains the curvarc-length parameterization, and when $\lambda\rightarrow+\infty$ the parmeterization tends to the arc-length parameterization.

%

%
%
%
	\vspace{-0.5cm}

\vspace{-1cm}
\section{Application to medical imaging : parameterization of bones}

\begin{figure}[!ht]
\centering
\includegraphics[height = 3cm]{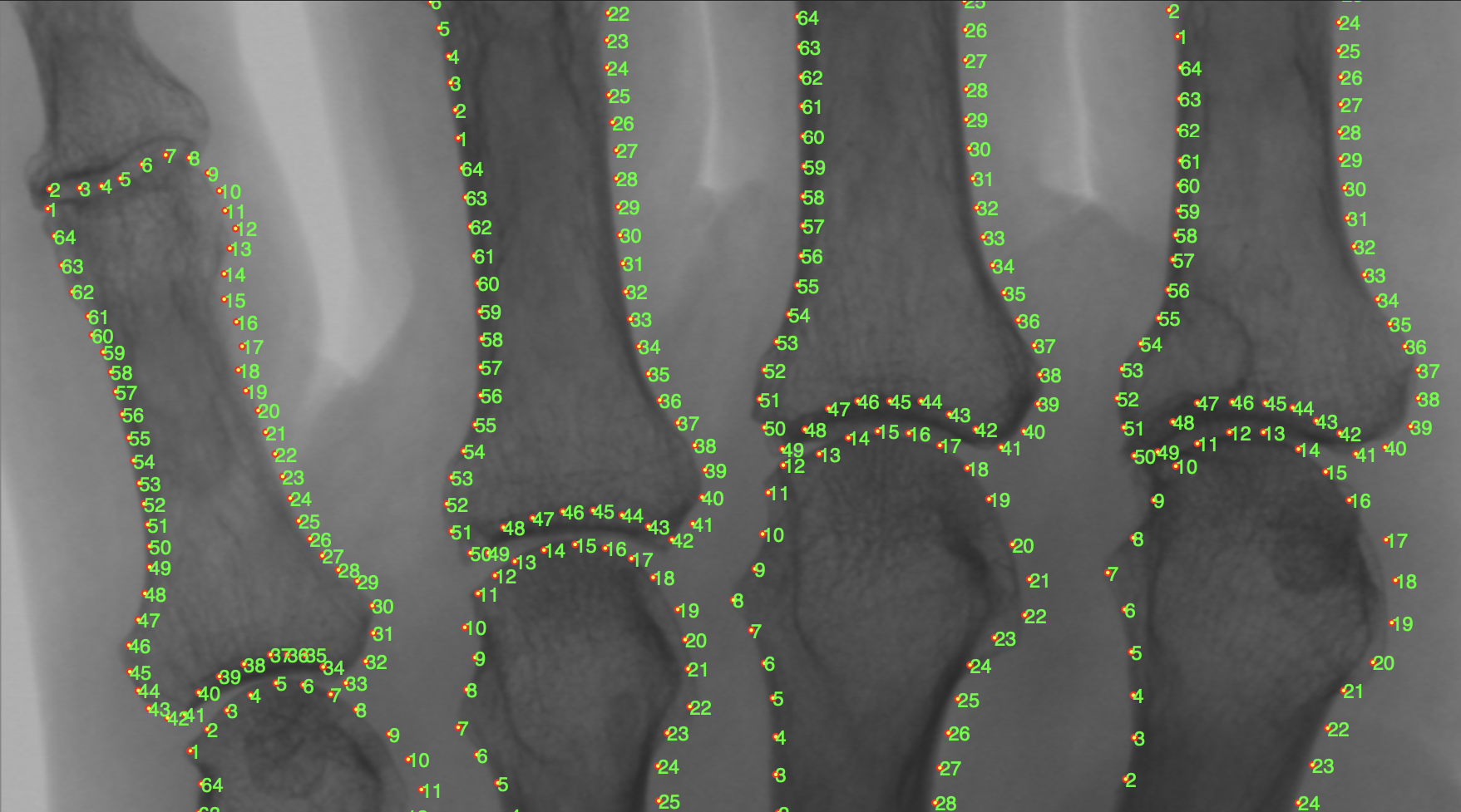}
\includegraphics[height = 3cm]{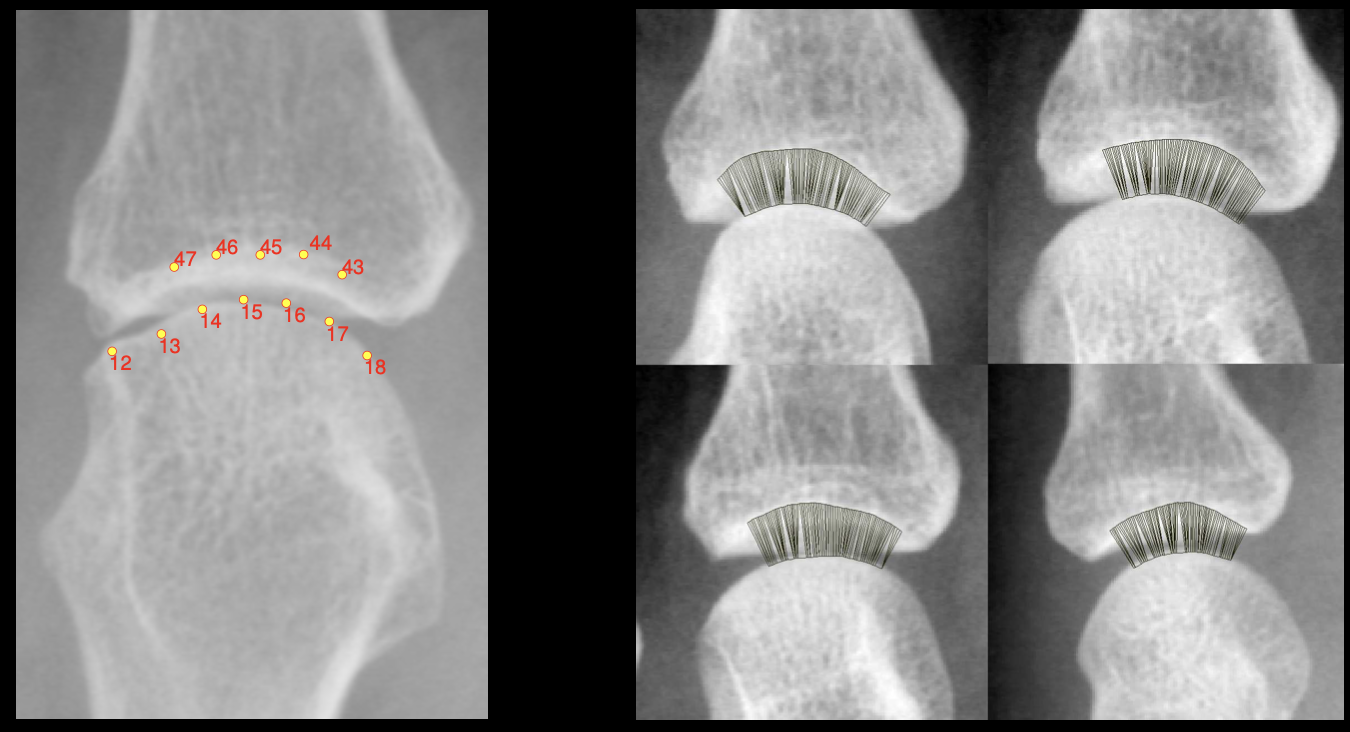}
\caption{Landmarks on bones used to measure joint space (courtesy of \cite{Langs})}
\label{landmarks}
\end{figure}
\begin{figure}[!ht]
 		\centering
		\includegraphics[height = 3cm]{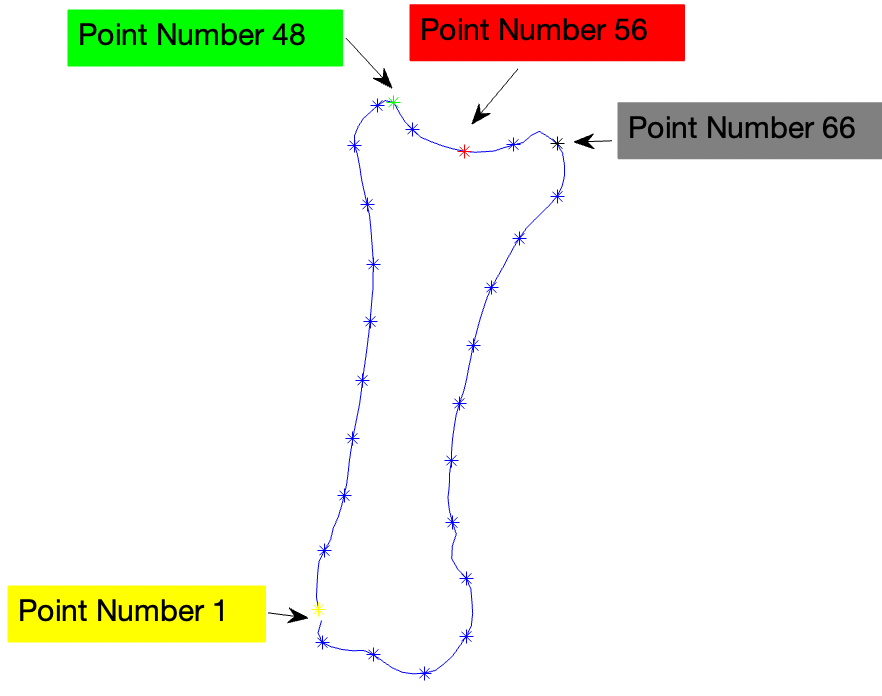}
		\includegraphics[height = 3cm]{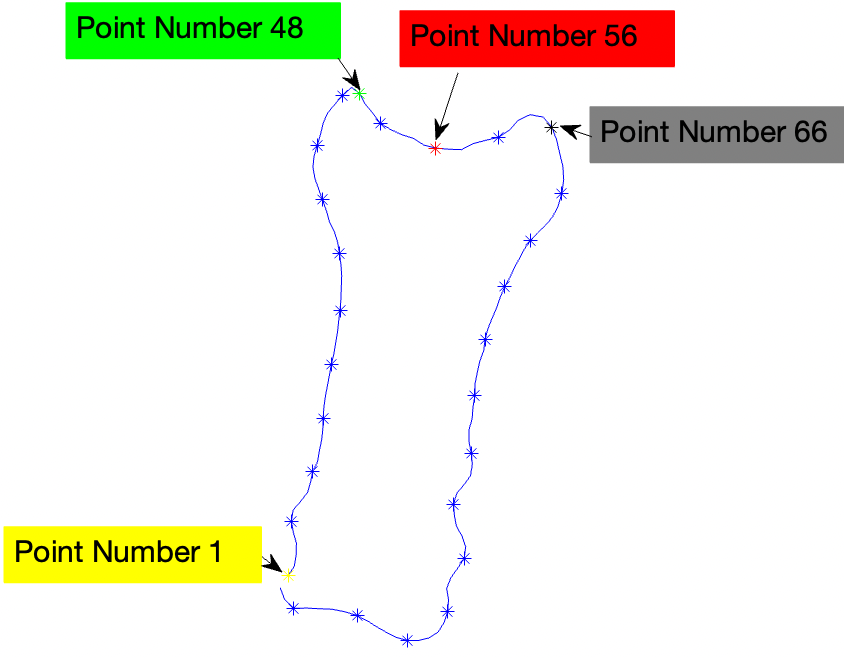}
		\includegraphics[height = 3cm]{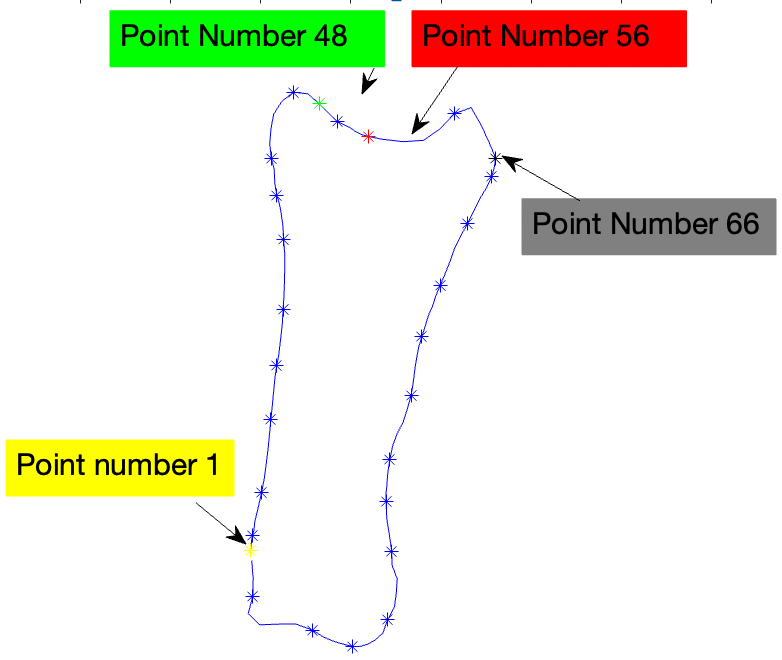}
		\caption{Point correspondance on 3 different bones using the method of \cite{Th}}
		\label{bones}
		\end{figure}

In the analysis of diseases like Rheumatoid Arthritis, one uses X-ray scans to evaluate how the disease affectes the bones. One effect of Rheumatoid Arthritis is erosion of bones, another is joint shrinking \cite{Langs}. In order to measure joint space, one has to solve a point correspondance problem. For this, one uses landmarks along the contours of bones as in Fig.\ref{landmarks}. 
These landmarks have to be placed at the same anatomical positions for every patient.  Below they are placed using  a method by Hans Henrik Thodberg \cite{Th}, based on minimum description length which minimizes the description of a PCA model capturing the variability of the landmark positions. 
For instance in Fig.~\ref{bones} left, the landmark number 56 should always be in the middle of the head of the bone because it is used to measure the width between two adjacent bones in order to detect rheumatoid arthritis.

%
Although the method by Hans Henrik Thodberg gives good results, it is computationally expensive. In this paper we propose to recover similar results with an quicker algorithm. It is based on the fact that any geometrically meaningful parameterization of a contour can be expressed using the arc-length measure and the curvature of the contour, which are the only geometric invariants of a 2D-curve (modulo translation and rotation). It follows that the parameterization calculated by Thodberg's algorithm should be recovered as a parameterization expressed using arc-length and curvature.
We investigate a  2 parameter family of parameterizations defined by 
\begin{equation}\label{u}
u(s) = \frac{\int_0^s (c* L + |\kappa(s)|^{\lambda}) ds}{\int_0^1 (c*L + |\kappa(s)|^\lambda) ds}
\end{equation}where $c$ and $\lambda$ are positive parameters and where $L$ is the length of the curve and $\kappa$ its curvature function.
We recover an analoguous parameterization to the one given by Thodberg's algorithm with $c = 1$ and $\lambda = 7$ at real-time speed (gain of 2 order of magnitude). 
 \begin{figure}[!ht]
 		\centering
\includegraphics[width = 10.5cm]{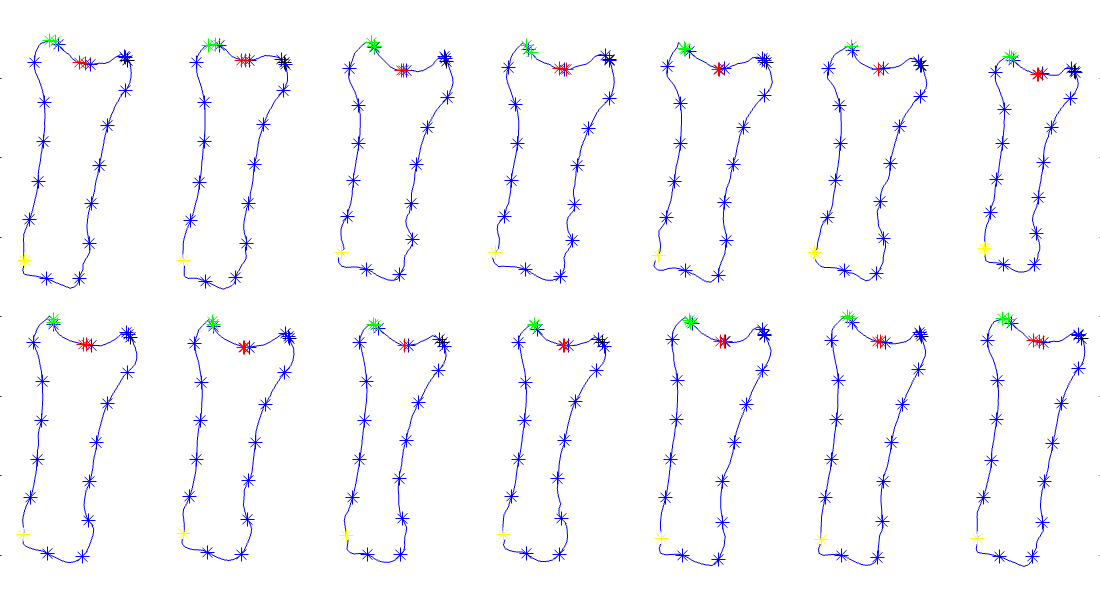}	
\caption{14 bones parameterized by Thodberg's algorithm on one hand and the parameterization defined by \eqref{u} with $c = 1$ and $\lambda = 7$ on the other hand (the two parameterizations are superposed). The colored points corresponds to points labelled $1$, $48$, $56$, $66$. They overlap for the two methods.}
\end{figure}


\section{Conclusion}

We proposed diverse canonical parameterization of 2D-contours, which are expressed using  arc-length and curvature of  curves. The curvature-length parameterization and the curvarc-length parameterization are very natural examples, since they corresponds to a constant-speed moving frame in $SO(2)$ and $SE(2)$.
We present an application to the problem of point correspondance in medical imaging consisting of labelling automatically keypoints along the contour of bones. We recover an analoguous parameterization to the one proposed by Thodberg at real-time speed. Having a two-parameter family of parameterizations at our disposal, a fine-tuning can be applied on top of our results in order to improve the point correspondance further.
%

%
%
%
%

\end{document}